\documentstyle [12pt,epsf, graphicx, amssymb]{article}

\begin{document}
\def\l{\lambda}
\def\m{\mu}
\def\a{\alpha}
\def\b{\beta}
\def\g{\gamma}
\def\d{\delta}
\def\e{\epsilon}
\def\o{\omega}
\def\O{\Omega}
\def\v{\varphi}
\def\t{\theta}
\def\r{\rho}
\def\bs{$\blacksquare$}
\def\bp{\begin{proposition}}
\def\ep{\end{proposition}}
\def\bt{\begin{th}}
\def\et{\end{th}}
\def\be{\begin{equation}}
\def\ee{\end{equation}}
\def\bl{\begin{lemma}}
\def\el{\end{lemma}}
\def\bc{\begin{corollary}}
\def\ec{\end{corollary}}
\def\pr{\noindent{\bf Proof: }}
\def\note{\noindent{\bf Note. }}
\def\bd{\begin{definition}}
\def\ed{\end{definition}}
\def\C{{\mathbb C}}
\def\P{{\mathbb P}}
\def\Z{{\mathbb Z}}
\def\d{{\rm d}}
\def\deg{{\rm deg\,}}
\def\deg{{\rm deg\,}}
\def\arg{{\rm arg\,}}
\def\min{{\rm min\,}}
\def\max{{\rm max\,}}

\newcommand{\norm}[1]{\left\Vert#1\right\Vert}
\newcommand{\abs}[1]{\left\vert#1\right\vert}

\newcommand{\set}[1]{\left\{#1\right\}}
\newcommand{\setb}[2]{ \left\{#1 \ \Big| \ #2 \right\} }

\newcommand{\IP}[1]{\left<#1\right>}
\newcommand{\Bracket}[1]{\left[#1\right]}
\newcommand{\Soger}[1]{\left(#1\right)}

\newcommand{\Integer}{\mathbb{Z}}
\newcommand{\Rational}{\mathbb{Q}}
\newcommand{\Real}{\mathbb{R}}
\newcommand{\Complex}{\mathbb{C}}

\newcommand{\eps}{\varepsilon}
\newcommand{\To}{\longrightarrow}
\newcommand{\varchi}{\raisebox{2pt}{$\chi$}}

\newcommand{\E}{\mathbf{E}}
\newcommand{\Var}{\mathrm{var}}

\def\squareforqed{\hbox{\rlap{$\sqcap$}$\sqcup$}}
\def\qed{\ifmmode\squareforqed\else{\unskip\nobreak\hfil
\penalty50\hskip1em\null\nobreak\hfil\squareforqed
\parfillskip=0pt\finalhyphendemerits=0\endgraf}\fi}

\renewcommand{\th}{^{\mathrm{th}}}
\newcommand{\Dif}{\mathrm{D_{if}}}
\newcommand{\Difp}{\mathrm{D^p_{if}}}
\newcommand{\GHF}{\mathrm{G_{HF}}}
\newcommand{\GHFP}{\mathrm{G^p_{HF}}}
\newcommand{\f}{\mathrm{f}}
\newcommand{\fgh}{\mathrm{f_{gh}}}
\newcommand{\T}{\mathrm{T}}
\newcommand{\K}{^\mathrm{K}}
\newcommand{\PghK}{\mathrm{P^K_{f_{gh}}}}
\newcommand{\Dig}{\mathrm{D_{ig}}}
\newcommand{\for}{\mathrm{for}}
\newcommand{\End}{\mathrm{end}}

\newtheorem{th}{Theorem}[section]
\newtheorem{lemma}{Lemma}[section]
\newtheorem{definition}{Definition}[section]
\newtheorem{corollary}{Corollary}[section]
\newtheorem{proposition}{Proposition}[section]

\begin{titlepage}

\begin{center}

\topskip 5mm

{\LARGE{\bf {Smooth rigidity and Remez inequalities  via Topology of level sets}}}

\vskip 8mm

{\large {\bf Y. Yomdin}}

\vspace{6 mm}
\end{center}

{The Weizmann Institute of Science, Rehovot 76100, Israel}

\vspace{2 mm}

{e-mail: yosef.yomdin@weizmann.ac.il}

\vspace{1 mm}

\vspace{1 mm}
\begin{center}

{ \bf Abstract}
\end{center}

{\small A smooth rigidity inequalitiy provides an explicit lower bound for the $(d+1)$-st derivatives of a smooth function $f$, which holds, if $f$ exhibits certain patterns, forbidden for polynomials of degree $d$. The main goal of the present paper is twofold: first, we provide an overview of some recent results and questions related to smooth rigidity, which recently were obtained in Singularity Theory, in Approximation Theory, and in Whitney smooth extensions. Second, we prove some new results, specifically, a new Remez-type inequality, and on this base we obtain a new rigidity inequality. In both parts of the paper we stress the topology of the level sets, as the input information. Here are the main new results of the paper:

\smallskip

Let $B^n$ be the unit $n$-dimensional ball. For a given integer $d$ let $Z\subset B^n$ be a smooth compact hypersurface with $N=(d-1)^n+1$ connected components $Z_j$. Let $\mu_j$ be the $n$-volume of the interior of $Z_j$, and put
$\mu=\min \mu_j, \ j=1,\ldots, N$. Then for each polynomial $P$ of degree $d$ on ${\mathbb R}^n$ we have
$$
\frac{\max_{B^n}|P|}{\max_{Z}|P|}\le (\frac{4n}{\mu})^d.
$$
As a consequence, we provide an explicit lower bound for the $(d+1)$-st derivatives of any smooth function $f$, which vanishes on $Z$, while being of order $1$ on $B^n$ (smooth rigidity)}:
$$
||f^{(d+1)}||\ge \frac{1}{(d+1)!}(\frac{4n}{\mu})^d.
$$
We also provide an interpretation, in terms of smooth rigidity, of one of the simplest versions of the results in \cite{Ler.Ste}.

\vspace{1 mm}
\end{titlepage}

\newpage


\section{Introduction}\label{Sec:Intro}
\setcounter{equation}{0}

Let $f(x)$ be a smooth function on the unit $n$-dimensional ball $B^n$.
A ``rigidity inequality'' for $f$ is an explicit lower bound for the $(d+1)$-st derivative of $f$, which holds, if $f$ exhibits certain patterns, forbidden for polynomials of degree $d$.

\medskip

We expect rigidity inequalities to be valid for those polynomial behavior patterns, which are stable with respect to smooth approximations. At present many such important patterns are known (see \cite{Ler.Ste,Yom1,Yom.Com,Yom2}). However, translation of the known ``near-polynomiality'' results into ``rigidity inequalities'' usually is not straightforward, and many new questions arise.

\medskip

Smooth rigidity  inequalities naturally form a certain domain of interrelated results and questions in Smooth Analysis. They are closely related to some other important mathematical fields, including Singularity Theory, Approximation Theory, Real Algebraic Geometry, and Whitney extension problems.

The main goal of the present paper is twofold: first, we provide an overview of some recent results and questions related to smooth rigidity. We start with Remez-type inequalities (\cite{Yom,Yom2,Yom3}), continue with recent exciting developments in the classical Whitney smooth extension theory (\cite{Bru.Shv,Fef,Fef.Kla}), and conclude with a recent important development in global Singularity Theory, achieved in \cite{Ler.Ste}.

\smallskip

Second, we prove some new results. This includes a new Remez-type inequality. On this base, via the results of \cite{Yom3}, we obtain a new Smooth rigidity inequality. We also want to illustrate in more detail some results of  \cite{Ler.Ste} and their consequences for smooth rigidity. For this purpose we give a direct proof of one very special case of the general results of \cite{Ler.Ste}. On this base we provide also the corresponding rigidity inequality. In all these new (or almost new) results the stress is on the topology of the level sets, as the input information.

\smallskip

The paper is organized as follows: Sections \ref{Sec:Rig.Remez} - \ref{Sec:Ler} form a review part of the paper. Specifically, in Sections \ref{Sec:Rig.Remez} and \ref{Sec:Remez.Const} we provide an overview of some results of \cite{Yom3}. In Section \ref{Sec:Rigid.Whitney} we shortly outline some connections of smooth rigidity with the recent important progress in the Whitney smooth extension problem (\cite{Bru.Shv,Fef,Fef.Kla}). In Section \ref{Sec:Ler} we provide a short overview of the results of \cite{Ler.Ste}.

\smallskip

Sections  \ref{Sec:Remez.type} and \ref{sec:rigidity.via.topology} provide new results: in Section \ref{Sec:Remez.type} we obtain a new Remez-type inequality, wich uses both the topological and measure information. On this base a corresponding rigidity inequality is obtained.  In Section \ref{sec:rigidity.via.topology} we give a proof of one result in the sprit of \cite{Ler.Ste}, and provide on this base the corresponding rigidity  inequality.

\section{Smooth rigidity - some background}\label{Sec:Rig.Remez}
\setcounter{equation}{0}

In this section we summarize, following \cite{Yom3}, some basic facts concerning smooth rigidity. One of possible specific setings of the Smooth rigidity problem was proposed in \cite{Yom3}. It assumes as an input data a certain closed subset $Z$ of the zeroes set $Y(f)$ of $f$. The following definition was given in \cite{Yom3}:

\smallskip

Let $f: B^n \rightarrow {\mathbb R}$ be a $d+1$ times continuously differentiable function on the unit closed ball $B^n \subset {\mathbb R}^n$. For $l=0,1,\ldots,d+1$ put
$$
M_l(f)=\max_{z\in B^n} \Vert f^{(l)}(z) \Vert,
$$
where the norm of the $l$-th derivative $f^{(l)}(z)$  of $f$ is defined as the sum of the absolute values of all the partial derivatives of $f$ of order $l$.

\medskip

For $Z\subset B^n$ let $W_d(Z)$ denote the set of $C^{d+1}$ smooth functions $f(z)$ on $B^n$, vanishing on $Z$, with $M_0(f)=1$.

\bd\label{def:rigidity}
For $Z\subset B^n$ we define the $d$-th rigidity constant ${\cal RG}_d(Z)$ as
$$
{\cal RG}_d(Z)=\inf_{f\in W_d(Z)}M_{d+1}(f).
$$
\ed
By this definition we get immediately $M_{d+1}(f)\ge {\cal RG}_d(Z)$ for any $f(z)$ on $B^n$, vanishing on $Z$, with $M_0(f)=1$. Our goal is to estimate ${\cal RG}_d(Z)$ in terms of accessible geometric features of $Z$. The rigidity constant
${\cal RG}_d(Z)$ is closely related to some of the central notions in the modern smooth extension theory (see \cite{Bru.Shv,Fef,Fef.Kla} and references therein). We give more comments on this connection in Section  \ref{Sec:Rigid.Whitney}.

\medskip

The following one-dimensional example illustrates some patterns of smooth rigidity. Start with a basic property of polynomials: a nonzero univariate polynomial $P(x)$ of degree $d$ can have at most $d$ real zeros. Here is  the corresponding rigidity result (well-known in various forms). It bounds from below the $(d+1)$-st derivative of a function $f$, which violates this property:

\bp\label{prop:d.zeroes}
For each $(d+1)$-smooth function $f(x)$ on $[-1,1]$, with $\max_{[-1,1]}|f(x)|\ge 1$ and with $d+1$ or more distinct zeroes on $[-1,1]$, we have
$$
\max_{[-1,1]}|f^{(d+1)}(x)|\ge \frac{(d+1)!}{2^{d+1}}.
$$
\ep
A short proof can be given using the Lagrange remainder formula in the polynomial interpolation of $f$ on its zeroes , or (more or less equivalently) via divided finite differences (see \cite{Yom3}).

\smallskip

In terms of the rigidity constant ${\cal RG}_d(Z)$ Proposition \ref{prop:d.zeroes} implies

\bp\label{prop:d.points}
For any $Z\subset B^1$ we have ${\cal RG}_d(Z)\ge \frac{(d+1)!}{2^{d+1}},$ if $Z$ consists of at least $d+1$ different points, and ${\cal RG}_d(Z)=0$ if $Z$ consists of at most $d$ different points.
\ep

The corresponding ``near-polynomiality'' result is the following:

\bc\label{cor:zeroes}
Any $(d+1)$-smooth function $f(x)$ on $[-1,1]$, with
$$
\max_{[-1,1]}|f(x)|\ge 1, \  \  \  \max_{[-1,1]}|f^{(d+1)}(x)|< \frac{(d+1)!}{2^{d+1}}
$$
has at most $d$ zeroes in $[-1,1]$.
\ec

\medskip

In higher dimensions the powerful one-dimension tools such as Lagrange's remainder formula, and divided finite differences, are not directly applicable. Still, Proposition \ref{prop:d.zeroes} implies, via line sections, the following

\bp\label{prop:Z.interior}
For any $Z \subset B^n$ with a non-empty interior,
$$
{\cal RG}_d(Z) \ge \frac{(d+1)!}{2^{d+1}}.
$$
\ep
\pr
Let $f\in U_d(Z)$. Fix a certain point $x_1$ with $|f(x_1)|=1$, fix $x_2$ in the interior of $Z$, and let $\ell$ be the straight line through $x_1,x_2$. The restriction $\bar f$ of $f$ to the intersection of $\ell$ with $B^n$ has an entire interval of zeroes near $x_2$, and it satisfies $M_0(\bar f)=1$. Applying Proposition \ref{prop:d.zeroes} to $\bar f$ completes the proof. $\square$

\smallskip

Considering restrictions of $f$ to the straight lines $\ell$ as above, and applying Proposition \ref{prop:d.zeroes}, we obtain in \cite{Yom1}, in particular, the following ``near-polynomiality'' result:

\bt\label{thm:poly.like.old}(\cite{Yom1})
Let $f(x)$ be a smooth function on $B^n$, with $M_0(f)=1$. If $M_d(f)\le 2^{-d-1}$, for some $d \ge 1,$ then the set of zeroes $Y$ of $f$ is contained in a countable union of smooth hypersurfaces, and the $(n - 1)$-Hausdorff measure of $Y$ is bounded by a constant, depending only on $n$ and $d$.
\et
In turn, we immediately obtain a certain multi-variate rigidity inequality: if for some $d \ge 1$  the set of zeroes $Z$ of a smooth function $f$ on $B^n$ of sup-norm $1$ violates the restrictions of Theorem \ref{thm:poly.like.old}, then $M_d(f)$ is at least $2^{-d-1}$.

\smallskip

However, the above approach to estimating ${\cal RG}_d(Z)$, based on restrictions of $f$ to some straight lines $\ell$, works only in situations, where ``many'' straight lines intersect $Z$ at ``many'' points. Essentially, (via the integral-geometric interpretation) this is the case of sets $Z\subset B^n$, containing hypersurfaces of a big Hausdorff $(n-1)$-measure. Still, for many natural classes of candidate zero sets $Z$, this condition is not satisfied. For example, this is the case for $Z$ being a finite set (unless, by a rare coincidence, many points of $Z$ lie on the same straight line). The same remains true if we replace each point of such a finite set $Z$ with a small hypersurface around it.

\medskip

The main goal of \cite{Yom3} was to develop a pure multi-dimensional approach to smooth rigidity, based on polynomial Remez-type inequalities (which compare the maxima of a polynomial on the closed unit ball, and on its closed subset $Z$). Loosely speaking, {\it one of the main results of \cite{Yom3} was that the $d$-rigidity of a set $Z$ is approximately the reciprocal $d$-Remez constant of $Z$}.

\section{Remez constant}\label{Sec:Remez.Const}
\setcounter{equation}{0}

Let us recall the definition and some basic properties of the Remez (or Lebesgue, or norming) constant. See \cite{Bru.Yom,Bru.Gan,Erd,Yom} for more details and references.

\bd\label{Remez.constant}
For a set $Z\subset B^n \subset {\mathbb R}^n$ and for each $d\in {\mathbb N}$ the Remez constant ${\cal R}_d(Z)$ is the minimal $K$
for which the inequality
$$
\sup_{B^n}\vert P \vert \leq K \sup_{Z}\vert P \vert
$$
is valid for any real polynomial $P(x)=P(x_1,\dots,x_n)$ of degree $d$.
\ed
Clearly, we always have ${\cal R}_d(Z)\ge 1.$ For some $Z$ the Remez constant ${\cal R}_d(Z)$ may be equal to $\infty$. In fact, ${\cal R}_d(Z)$ is infinite if and only if $Z$ is contained in the set of zeroes
$$
Y_P=\{x\in {\mathbb R}^n, \ | \ P(x)=0\}
$$
of a certain polynomial $P$ of degree $d$. We call sets $Z$ with finite ${\cal R}_d(Z)$ $d$-norming. We use also the reciprocal Remez constant $\hat {\cal R}_d(Z):=\frac{1}{{\cal R}_d(Z)}.$

\subsection{Rigidity constant via Remez constant}\label{Sec:rigid.via.topol}

An important initial observation, connecting the Remez and rigidity constants is:

\bl\label{lem:R.is.Inf.11}(\cite{Yom3})
${\cal RG}_d(Z)=0$ if and only if $\hat {\cal R}_d(Z)=0.$
\el
See \cite{Yom3} for the proof. The following is one of the main results of \cite{Yom3}. It is based, in particular, on \cite{Yom2}:

\bt\label{thm:main11}(\cite{Yom3})
For any $Z \subset B^n$, \ \ \ we have \ \  $\frac{(d+1)!}{2}\hat {\cal R}_d(Z)\le {\cal RG}_d(Z)$.
\et
This lower bound is valid for any $Z$, and it is sharp, up to constants (depending only on $n$ and $d$, and on the separation between the point), for finite sets, as Theorem \ref{thm:main2} below shows. However, we cannot expect {\it an upper bound} of the form

\be\label{eq:both.sides1}
{\cal RG}_d(Z)\le C(n,d) \hat {\cal R}_d(Z),
\ee
for some constant $C(n,d)$ depending only on $n$ and $d$, to be valid in general: indeed, by Proposition \ref{prop:Z.interior}, for any $Z \subset B^n$ with a non-empty interior, ${\cal RG}_d(Z) \ge \frac{(d+1)!}{2^{d+1}}.$ On the other hand, sets $Z$ with a non-empty interior may have arbitrarily small Remez constant $\hat {\cal R}_d(Z)$. For example, let $P$ be a polynomial of degree $d$ with $M_0(P)=1$, and  for some $\eta>0$ let $Z$ be the $\eta$-sublevel set of $P$, i.e. $Z=\{z\in B^n, \ |P(z)|\le \eta\}$. Clearly, $Z$ has a non-empty interior, and we have $\hat {\cal R}_d(Z) \le \eta$.

\medskip

Still, for some important types of sets $Z$ an upper bound for the rigidity through the Remez constant holds. In \cite{Yom3} we prove it for finite sets $Z$:

\bt\label{thm:main2}(\cite{Yom3})
Let $Z \subset B^n$ be a finite set, and let $\rho$ be the minimal distance between the points of $Z$. Then
$$
\frac{(d+1)!}{2}\hat {\cal R}_d(Z)\le {\cal RG}_d(Z)\le \frac{C(n,d)}{\rho^{d+1}}\hat {\cal R}_d(Z).
$$
\et
This theorem can be considered as a generalization of Proposition \ref{prop:d.points} to higher dimensions.

\smallskip

In dimensions $2$ and higher we have finite sets $Z$ with positive but arbitrarily small $\hat{\cal R}_d(Z),$ and with $\rho$, uniformly bounded from below. For such sets the upper bound of Theorem \ref{thm:main2} is meaningful. One of the simplest examples is a plane triangle $Z_h$, defined as
$$
Z_h=\{(-\frac{1}{2},0),(0,h),(\frac{1}{2},0)\}.
$$
Easy computation shows that $\hat {\cal R}_1(Z_h)=\frac{h}{2}$.

\section{Smooth rigidity and Whitney extensions}\label{Sec:Rigid.Whitney}
\setcounter{equation}{0}

In this section we discuss, quite informally, some very important connections between the smooth rigidity and Whitney smooth extension (see \cite{Whi1,Whi2,Whi3,Bru.Shv,Fef,Fef.Kla}). In fact, our Definition \ref{def:rigidity} of the $d$-rigidity ${\cal RG}_d(Z)$ is a special case of one of the main notions in the Whitney smooth extension theory. Indeed, in the Whitney $C^m$-extension problem we consider a closed subset $Z\subset B^n$, and a function $\tilde f$ on $Z$. The question is whether $\tilde f$ is extendable to a $C^m$-smooth function $f$ on $B^n$, and if so,  what is the minimal $C^m$-norm of the smooth extensions $f$ of $\tilde f$ to $B^n$. In our Definition \ref{def:rigidity} we just assume that $\tilde f \equiv 0$ on $Z$, add the requirement $M_0(f)=1$, and ask for the minimal $M_d(f)$ of the extensions $f$.

\smallskip

In dimension $n=1$ the Whitney extension theorem of \cite{Whi2} provides the complete and explicit answer to the $C^m$ - extension question: it is possible if and only if all the divided finite differences of $\tilde f$ on the subsets of $Z$ of cardinality at most $m+1$ are uniformly bounded. The minimal $C^m$ - norm of the extensions is also estimated through these finite differences.

\smallskip

As a consequence we can produce an explicit expression for the $d$-th rigidity constant ${\cal RG}_d(Z)$. Consider all the subsets $\tilde Z = \{z_0,z_1,\ldots,z_{d+1}\}$, with $z_0\in [-1,1]\setminus Z$ and $z_1,\ldots,z_{d+1}\in Z,$ and let $\Delta_{d+1}(\tilde Z)$ denote the $d+1$-st divided finite difference on $\tilde Z$ of the function $y_0=1, y_1=\ldots = y_{d+1}=0$. {\it Then the $d$-th rigidity constant ${\cal RG}_d(Z)$ can be estimated as the infinum over all the subsets $\tilde Z$ as above of $\Delta_{d+1}(\tilde Z)$}.
\smallskip

There is a fundamental difficulty in extending one-dimensional polynomial interpolation and smooth extension results to higher dimensions. This difficulty manifests itself in many ways, but for our purposes it can be shortly summarized as follows: {\it in dimensions greater than one there are no canonical divided finite differences}. Even the following most basic question, directly suggested by Whitney's one-dimensional results, was open for many years:

\smallskip

{\it In order to check a $C^m$-extensibility of $\tilde f$ on $Z\subset B^n$, and to estimate a minimal $C^m$-norm of the extension, is it enough to check only subsets of $Z$ of a fixed cardinality $N=N(n,m)$?}

\smallskip

A remarkable progress was achieved in the multi-dimensional  Whitney extension problem in the last two decades (see \cite{Bru.Shv,Fef,Fef.Kla} and references therein). In particular, the above question was ultimately answered positively in \cite{Fef}. As a result, we can provide an explicit expression for the $d$-th rigidity constant ${\cal RG}_d(Z)$ through the rigidity constants ${\cal RG}_d(\bar Z)$, where $\bar Z$ runs over all the finite subsets of $Z$ of cardinality $N=N(n,d)$. However, in a strict contrast with the one-dimensional case,  the rigidity constants ${\cal RG}_d(\bar Z)$ for finite sets $\bar Z$ of cardinality $N$ are not given anymore by a simple formula.

\smallskip

Accordingly, we come to an important problem of providing some ``computable'' bounds for ${\cal RG}_d(\bar Z)$ for {\it finite sets} $\bar Z$ in terms of their geometry. Some results in this direction were obtained in \cite{Yom3}. In particular, besides the results in Sections \ref{Sec:Rig.Remez} and \ref{Sec:Remez.Const} above, we have the following ``density'' lower bound for the rigidity ${\cal RG}_d(Z)$:

\bt\label{thm:main3} (\cite{Yom3})
Let $Z \subset B^n$ be a finite set, and let $\rho$ be the minimal distance between the points of $Z$. Assume that the cardinality $M=|Z|$ satisfies
$M > (4d)^n(\frac{1}{\rho})^{n-1}$. Then
$$
0 < \frac{(d+1)!}{2}\left ( \frac{M\rho^n- (4d)^n\rho}{4n}\right )^d \le {\cal RG}_d(Z).
$$
\et
Informally, if in resolution $\rho$ the set $Z$ ``looks more massive than an algebraic hypersurface of degree $d$'', then its $d$-th rigidity constant is positive.

\smallskip

Some additional results, in terms of the asymptotic behavior of the covering number of $\bar Z$, were obtained in \cite{Yom4}. We plan to rigorously present our results on connections between Whitney extension and smooth rigidity, shortly outlines in this section, separately.

\section{What is the degree of a smooth hypersurface (after A. Lerario and M. Stecconi, \cite{Ler.Ste})}\label{Sec:Ler}
\setcounter{equation}{0}

In this section we continue the review part of the preset paper. We discuss, from the point of view of smooth rigidity, an important recent progress in Global Singularity Theory, achieved in \cite{Ler.Ste}. Informally, the results of \cite{Ler.Ste} claim that, under appropriate transversality conditions, singular loci of smooth functions are isotopic to the corresponding singular loci of algebraic polynomials of explicitly bounded degree. It is shown, in particular, that if $f$ as above vanishes {\it transversally} on a smooth manifold $Z\subset B^n$, and if high-order derivatives of $f$ are small enough, then $Z$ is smoothly isotopic to a component of a zero set of a polynomial $P$ of a ``small'' degree. We expect that the results of \cite{Ler.Ste} will lead to some important ``rigidity via topology'' inequalities for  singular loci of smooth functions.

\smallskip

Let's quote \cite{Ler.Ste} and state one of its main results more accurately. For a smooth $f: B^n\to {\mathbb R}$ with a compact and smooth zero set $Y(f)$ the following problem is considered: what is the smallest degree of a polynomial $P$ whose zero set $Y(P)$ is diffeomorphic to $Y(f)$? More generally, for $W$ a closed semialgebraic subset of the jet space $J^r(B^n,{\mathbb R}^q)$ we consider the (type $W$) singularity $Y_W(f)$, defined as
$$
Y_W(f)=(j^r f)^{-1}(W),
$$
where $j^r f: B^n\to J^r(B^n,{\mathbb R}^q)$ is the $r$-jet extension of $f$. Certainly, some transversality assumptions are necessary in order to study the topology of the singular sets $Y_W(f)$. Otherwise, even in the simplest case of the zero level hypersurfaces of $f$ any closed subset of $B^n$ may appear as $Y(f)$ for a certain $C^\infty$ function $f$.

\smallskip

Define the discriminant set $\Delta_W$ as the set of all $f\in C^{r+1},$ whose $r$-jet extension $j^r f: B^n\to J^r(B^n,{\mathbb R}^q)$ is not transversal to $W$. Quantitatively, the transversality of $f$ to $W$ is measured by the distance $\delta_W(f)$ from $f$ to $\Delta_W$ in the space of $C^{r+1}$-smooth functions.

\bt\label{thm:Ler1}(\cite{Ler.Ste}, Theorem 1)
Let $W\subset J^r(B^n,{\mathbb R}^q)$ be closed semialgebraic. For every $f\in C^{r+2}(B^n,{\mathbb R}^q)$ with $j^rf$ transversal to $W$ there exists a polynomial map $P=(P_1,\ldots,P_q)$ with
$$
\deg (P_i) \le C(r) \max \left \{ r+1, \frac{||f||_{C^{r+2}}}{\delta_W(f)} \right \},
$$
such that the (type $W$) singularity $Y_W(f)$ is smoothly isotopic to $Y_W(P)$.
\et
This result allows the authors to immediately bound the topological complexity of $Y_W(f)$ (see \cite{Ler.Ste}, Corollary 2).

\smallskip

From the smooth rigidity poin of view, Theorem \ref{thm:Ler1} immediately implies the following result, where we put $\bar  C(r) = \frac {1}{ C(r)}$:

\bt\label{thm:Ler.Rigid.1}
For $f$ as above, and for a certain closed semialgebraic $W\subset J^r(B^n,{\mathbb R}^q)$, assume that the (type $W$) singularity $Y_W(f)$ is not smoothly isotopic to $Y_W(P)$ for any polynomial map $P$ of a given degree $d\ge r+1$. Then
$$
||f||_{C^{r+2}}\ge \bar  C(r) \cdot d \cdot \delta_W(f).
$$
\et
Of course, the condition of $Y_W(f)$ not being smoothly isotopic to $Y_W(P)$ for any polynomial map $P$ of degree $d\ge r+1$ can be replaced by stronger conditions in terms of topological complexity of $Y_W(P)$ (for example, in terms of the sum of its Betti numbers).

\smallskip

If we take the degree $d$ as the ``measure of complexity'' of polynomials, then Theorem \ref{thm:Ler.Rigid.1} implies the {\it linear grows of the norm $||f||_{C^{r+2}}$ with $d$}, at least. In would be interesting to make this very initial observation more rigorous.

\smallskip

We believe that an explicit and detailed study, in the above directions, of the {\it specific} singularity types $W$ is well justified. As an example, in Section \ref{sec:rigidity.via.topology} below we provide a strongly simplified (and restricted only to the zero hypersurfaces $Y(f)$) version of Theorem \ref{thm:Ler1}. Our goals there are to get explicit bounds, to estimate separately the $r+1$-st derivative of $f$, and to provide simple geometric arguments in the proof of this specific case.

\section{``Topological'' Remez-type inequality}\label{Sec:Remez.type}
\setcounter{equation}{0}

Starting with this section, we present the new results of the present paper. Remez-type inequalities provide an upper bound for ${\cal R}_d(Z)$ in terms of various ``computable'' characteristics of $Z$. In particular, the classical multi-dimensional Remez inequality (\cite{Bru.Gan}, \cite{Rem}, see also \cite{Erd}) uses the Lebesgue measure of $Z$. It reads as follows:

\bt\label{Remez.1}
For any measurable $Z\subset B^n$ we have
\be\label{Remez.ineq.n}
{\cal R}_d(Z) \ \leq \ T_d ({{1+(1-\lambda)^{1\over n}}\over {1-(1-\lambda)^{1\over n}}})\le (\frac{4n}{\lambda})^d.
\ee
Here $T_d(t)=cos(d \ arccos \ t)$ is the Chebyshev polynomial of degree $d$, and \ \ \ $\lambda= {{m_n(Z)}\over {m_n(B^n)}},$ with $m_n$ being the Lebesgue measure on ${\mathbb R}^n$. This inequality is sharp and for $n=1$ it coincides with the classical Remez inequality of \cite{Rem}.
\et
Some other examples of Remez-type inequalities, and a more detailed discussion can be found in \cite{Bru.Yom,Erd,Yom}. In particular, it was shown in \cite{Yom} that the Lebesgue measure can be replaced in Theorem \ref{Remez.1} with a more sensible geometric invariant $\o_{n,d}(Z)$, defined in terms of the covering numbers of $Z$. The invariant $\o_{n,d}(Z)$ always satisfies $m_n(Z)\le \o_{n,d}(Z)$, so its substitution instead of the Lebesgue measure into  (\ref {Remez.ineq.n}) can only improve the result. The invariant $\o_{n,d}(Z)$ allows us to distinguish between various discrete and even finite sets of different geometry and density.

\smallskip

However, there are natural classes of sets $Z$, for which Theorem \ref{Remez.1}, as well as its strengthening, where $m_n(Z)$ is replaced with $\o_{n,d}(Z)$, do not work. Consider, for instance, smooth compact hypersurfaces $Z\subset B^n$. Their $n$-measure is zero, and if their $n-1$-area is small, then also $\o_{n,d}(Z)=0$. Still $Z$ may be a $d$-norming set by ``topological'' reasons. Theorem \ref{thm:remez.topology1} below is one of the main results of the present paper. It provides a Remez-type inequality for smooth compact hypersurfaces $Z\subset B^n$, in terms of the number of their connected components, and of the $n$-volume of the {\it interiors of these components} (and not of the components themselves, as in Theorem \ref{Remez.1}).

\medskip

We state this result in a slightly more general way: let $U_j, \ j=1,\ldots,N,$ be nonintersecting compact connected domains in $B^n$ with nonempty interiors, and let $\mu_j=m_n(U_j)$ be the $n$-measure of $U_j$. It is convenient to assume that $\mu_1\ge \mu_2 \ge \ldots \ge \mu_N$. Put $Z_j=\partial U_j,$ to be the boundary of $U_j$, and let $Z=\cup_{j=1}^N Z_j$.

\medskip

Define $\bar d$ by $(\bar d-1)^n+1 \le N  < \bar d^n+1,$ and for each natural $d\le \bar d$ put $j_d=(d-1)^n+1$.

\bt\label{thm:remez.topology1}
For each $d \le \bar d$ we have
$$
{\cal R}_d(Z)\le (\frac{4n}{\mu_{j_d}})^d.
$$
\et
\pr


Let $d\le \bar d$ be fixed, and let $P$ be a polynomial of degree $d$ with $M_0(P)=1$. To prove Theorem \ref{thm:remez.topology1} it is sufficient to show that $\max_Z |P|\ge \kappa_d: = (\frac{\mu_{j_d}}{4n})^d.$ Thus we assume, in contrary, that $\max_Z |P| < \kappa_d,$ and bring this assumption to a contradiction.

\smallskip

We have the following lemma:

\bl\label{lem:max.Uj}
For each $U_j, \ j=1,\ldots,N,$ we have
$$
\max_{U_j}|P(x)|\ge ({\mu_j\over {4n}})^d.
$$
\el
\pr
By the classical multidimensional Remez inequality, given by Theorem \ref{Remez.1} above, and applied to the set $U_j$, we have
$$
1=M_0(P)\leq ({4n\over \ \mu_j})^d \max_{U_j}|P(x)|, \ \ \ \ \ or \ \ \ \ \ \max_{U_j}|P(x)|\ge ({\mu_j\over {4n}})^d.
$$
This completes the proof of Lemma \ref{lem:max.Uj}. $\square$

\medskip

Next we notice that if $\max_Z |P| < \kappa_d,$ then for each $j=1,\ldots, j_d$, the polynomial $P$ has a local maximum (or minimum) at a certain point $\bar x_j$ in the interior of $U_j$. Indeed, in this case, by Lemma \ref{lem:max.Uj}, the maximum of $|P|$ inside $U_j$ satisfies
$$
\max_{U_j}|P(x)|\ge ({\mu_j\over {4n}})^d \ge ({\mu_{j_d}\over {4n}})^d=\kappa_d,
$$
while $\max_Z |P| < \kappa_d$. Hence the maximum of $|P|$ in $U_j$ is strictly greater than its maximum on the boundary $Z_j$. In particular, the point $\bar x_j,$ where this maximum is achieved, is a critical point of $P$, i.e. $grad P(\bar x_j)=0$.

\medskip

Consequently, if $\max_{Z}|P(x)|< \kappa_d,$ then in the interior of each domain $U_j, \ j=1,\ldots,j_d,$ there is a critical point of $P$, which is a local maximum or a local minimum of $P$. Performing a small perturbation of $P$ we can assume that all the critical points of $P$ are non-degenerate.

\smallskip

It remains to bound from above the maximal possible total number of non-degenerated maxima and minima of a polynomial $P$ of degree $d$ of $n$ variables. By Bezout theorem, the total number of non-degenerated critical points of $P$, i.e. of the solutions of the system $\frac{\partial f}{\partial x_i}=0, \ i=1,\ldots,n$, cannot exceed $(d-1)^n$. This contradiction completes the proof of Theorem \ref{thm:remez.topology1}. $\square$

\medskip

Now, let $Z$ be as in Theorem \ref{thm:remez.topology1}. Combining  Theorems \ref{thm:remez.topology1} and \ref{thm:main11}, we obtain a ``topological'' rigidity inequality:

\bt\label{thm:rigidity}
For each $d \le \bar d$ we have
$$
{\cal RG}_d(Z)\ge \frac{(d+1)!}{2}(\frac{\mu_{j_d}}{4n})^d.
$$
\et


\subsection{Some examples and remarks}\label{Sec.Examples}

An immediate corollary of Theorem \ref{thm:remez.topology1} is the following:

\bc\label{cor:curves1}
Let the degree $d$ be given. Then for $N=(d-1)^n+1$ each set $Z=\cup_{j=1}^N Z_j$ as above is $d$-norming, and
$$
{\cal R}_d(Z)\le (\frac{4n}{\mu_{N}})^d.
$$
\ec
In a special case of  exactly one domain we have the following corollary of Theorem \ref{thm:remez.topology1}:

\bc\label{cor:curves}
Let $Z$ be the boundary of the compact connected domain $U$ in $B^n, \ n\ge 1$, with the $n$-volume of $U$ being $\mu$. Then we have
$$
{\cal R}_1(Z)\le \frac{4n}{\mu}.
$$
\ec

\bc\label{cor:curves2}
Let $Z_1,Z_2$ be the boundaries of the compact connected nonintersecting domains $U_1,U_2$ in $B^n, \ n\ge 1$, with the $n$-volumes of $U_1,U_2$ being $\mu_1 \ge \mu_2$. Then for $Z=Z_1\cup Z_2$ we have
$$
 {\cal R}_2(Z)\le (\frac{4n}{\mu_{2}})^2.
$$
\ec
\pr
For any $n\ge 1$ we have, by definition, $j_1=1, \ j_2=2$. Hence the result follows directly from Theorem \ref{thm:remez.topology1}. $\square$

\medskip

Theorem \ref{thm:remez.topology1} is sharp up to constants (depending on $n$ and $d$) with respect to the volume $\mu_{j_d}$ and with respect to the required number $N$ of the domains $U_j$. We give here in detail only the simplest example, for $n=2$ and $d=2$, i.e. we consider Corollary \ref{cor:curves2} instead of Theorem \ref{thm:remez.topology1}. Then we shortly discuss also the general case.

\smallskip

For a given $h>0$ consider a polynomial $P_h(x,y)=h^2x^2+y^2-\frac{1}{4}h^2$. The zero set $Y_h$ of $P_h(x,y)$ is the ellipse centered at the origin, with the semiaxes $\frac{1}{2}$ and
$\frac{h}{2}$ in the directions $Ox$ amd $Oy$, respectively. Now in Corollary \ref{cor:curves2} we put $U_1$ to be the interior of the ellipse $Y_h$, and $U_2$ to be the interior of the rectangle
$$
Q_h=\{(x,y)  \  \  |  \  \ - \frac{1}{4} \le x \le \frac{1}{4}, \  \  \  \frac{2h}{3} \le y \le \frac{3h}{4} \}.
$$
We see immediatly that just one oval in $Z$ is not enough to make $Z$ a $2$-norming set. Indeed, $Z_1$ is the zero set the polynomial $P_h$ of degree $2$. Next, the maximum of $|P_h|$ on $Q_h$ does not exceed $h^2$ while the maximum of $|P_h|$ on $B^2$ is at least $1- \frac{1}{4}h^2$. We conclude that for $Z=Z_1\cup Z_2$ we have
$$
 {\cal R}_2(Z)\ge \frac{1- \frac{1}{4}h^2}{h^2} \ge \frac{1}{2h^2}
$$
for small $h$. On the other hand, the smallest area of $U_1,U_2$ is of $U_2$, and we have $\mu_2=\frac{h}{48}$. Thus the bound of Corollary \ref{cor:curves2}, for $n=2$, takes the form
$$
 {\cal R}_2(Z)\le (\frac{4n}{\mu_{2}})^2=(\frac{8\cdot 48}{h})^2=\frac{147456}{h^2}.
$$
Therefore,  the power, with which the volume $\mu_{2}$ enters the bound, is accurate, while the bound itself is sharp, up to a constant.

\smallskip

The requirement of Theorem \ref{thm:remez.topology1} to have at least $N=(d-1)^n+1$ different disjoint domains $U_j$ remains relatively sharp also in higher degrees and dimensions.

\smallskip

Gonsider the following polynomial $P(x_1,\ldots,x_n)$ of degree $nd$:
$$
P(x_1,\ldots,x_n)=\prod_{i=1}^n Q(x_i),  \   \   \  Q(t)=(t-\eta_1)(t-\eta_2)\cdot \ldots \cdot (t-\eta_d),
$$
with $\eta_1,\cdots,\eta_d$ pairwise distinct numbers in the open interval $(-\frac{1}{\sqrt n}, \frac{1}{\sqrt n})$. The zero set $Y$ of the polynomial $P$ is the union of all the shifted coordinate hyperplanes ${x_i=\eta_j}, \  i=1,\ldots, n, \  j=1,\ldots, d$. On each  connected component of the complement $W=B^n \setminus Y$ the polynomial $P$ preserves its sign, and it changes its sign as the argument crosses $Y$.  The complement $W$ contains, in particular, $(d-1)^n$ adjoint cubes, on a half of them $P$ being positive, and on a half negative. Now fix a small positive number $\zeta$, which is a regular value of $P$, put $\bar P = P - \zeta$, and consider the sublevel set $V=\{x\in B^n, \  \bar P(x)\ge 0\}$.
We take as $U_j, \ j=1,\ldots,\frac{1}{2}(d-1)^n$ all the compact connected components of $V$, inside the cubes, were $P$ was positive. The boundaries $Z_j$ of $U_j$ are smooth compact hypersurfaces, contained in the zero set $\bar Y$ of $\bar P$. Therefore for $Z=\cup Z_j$ we have ${\cal R}_{nd}(Z)=\infty$. Replacing $d$ by $\frac {d}{n}$ we produce an example of $\hat N = \frac{1}{2}(\frac{d}{n}-1)^n$ disjoint connected domains $\hat U_j$ for which $\hat Z = \cup \hat Z_j$ is not a $d$-norming set. We conclude that the required number $N=(d-1)^n+1$ of the domains $U_j$ in Theorem \ref{thm:remez.topology1} is sharp in the degree $d$, up to a constant depending only on the dimension $n$.

\smallskip

We expect than the power $d$, with which the volume $\mu_{j_d}$ enters the bound, is accurate. However, the above construction immediately produces examples with $ {\cal R}_d(Z) \sim (\frac{1}{\mu_{j_d}})^{\frac{d}{n}}$ only. Indeed, we can consider the product of the polynomials $Q(x_i)$ as above, with the roots of the first one down-scaled to the size $h$. Then all the domains $U_j$ constructed belong to the strip $|x_1|\le h$. Repeating verbally the construction of the example after
Corollary \ref{cor:curves2}, we obtain ${\cal R}_d(Z) \sim (\frac{1}{h})^{\frac{d}{n}}$, while $\mu_{j_d}\sim h$.

\smallskip


\smallskip

Notice that for $d=1,2$ and for any $n$, the bound $(d-1)^n+1 = 1,2,$ on the number of the connected components of $Z$, is sharp. For $d\ge 3$, in many cases this bound can be improved, since by topological reasons, other critical points, beyond maxima and minima, must appear. Still, we are aware only of some partial estimates of the possible number of minima and maxima of real polynomials. In particular, one can show that the number of extrema of $P$ does not exceed $\frac{1}{2} d^n + O(d^{n-1})$. For n=3, in \cite{Bih}  there are examples with at least $\frac{13}{36} d^n$ extrema. Also the lower bound $(\frac{2^{n-1}}{n!}) d^n$ is known. (The author thanks E. Shustin for providing some references).

\smallskip

We consider an accurate estimate of the required number of the connected components of $Z$ in Theorem \ref{thm:remez.topology1} as an interesting question in real algebraic geometry, closely related not only to bounding the number of the extrema of $P$, but also to the topology and mutual position of its ovals.

\smallskip


\smallskip

Returning to the product polynomial $P(x_1,\ldots,x_n)=\prod_{i=1}^n Q(x_i)$, constructed above, we notice that $P$ can be considered as a ``poly-degree $d$'' polynomial, with respect to the appropriate Newton diagram. Many results of real algebraic geometry can be extended to such polynomials, as well as some Remez-type inequalities. We expect that the connection between these two topics, provided by Theorem \ref{thm:remez.topology1}, remain valid also for polynomials with a prescribed Newton diagram (and not only of a given degree).


\section{Topology of transversal level sets}\label{sec:rigidity.via.topology}
\setcounter{equation}{0}

In this section we illustrate in more detail the results of  \cite{Ler.Ste} and their consequences for smooth rigidity. For this purpose, and in order to provide as explicit geometric arguments as possible, we give below  a direct proof of one very special case of the general results of \cite{Ler.Ste}. On this base we provide also the corresponding rigidity inequality.

\smallskip

Let $f$ be a $C^{d+1}$-smooth function on $B^n,$ with $M_0(f)=1$. For a given $\gamma > 0$, a real number $c$ is called a $\gamma$-regular value of $f$, if for each $x\in B^n$ with $f(x)=c$, we have $||grad \ f(x)|| \ge \ \gamma$. In this case the level set
$$
Y_c(f)=\{x\in B^n, \ f(x)=c\}
$$
is a regular compact manifold of dimension $n-1$.

\smallskip

To simplify the presentation, and to avoid the boundary effects, we assume that $f$ does not vanish out of  the concentric ball  $B^n_{\frac{1}{2}}.$ Next we assume that for a certain $\gamma, \ 0<\gamma < 1$, zero is a $\gamma$-regular value of $f$. In the notations of \cite{Ler.Ste}, for $W$ consisting of the jets with the zero value, $\gamma$ is, essentially, the distance $\delta_W(f)$ of $f$ to the discriminant set $\Delta_W$ .

\smallskip

Thus the zero level set $Y_0(f)=\{x\in B^n, \ f(x)=0\}$ is  a compact smooth hypersurface in $B^n_{\frac{1}{2}}.$ Let $V_i, \ i=1,\ldots, q,$ denote the connected components of $Y_0(f)$.

\smallskip

We put $T=T(\gamma)=\min \{1,\frac{d!\gamma^2}{4C_3}\}$, where the constant $C_3=C_3(n,d)$ is defined below. Finally, let $P=P_d(f)$ be the Taylor polynomial at the origin of degree $d$ of $f$, and let $Y_0(P)$ be its zero set. The following result is (essentially) a very special case of Theorem 1.1 of \cite{Ler.Ste} (see also Section \ref{Sec:Ler} above):

\bt\label{thm:rigidity.topology}
If $M_{d+1}(f)\le T$ then the smooth hypersurface $Y_0(f)=\cup_{i=1}^q V_i $ is smoothly isotopic to a certain union $W=\cup_{i=1}^q W_i $ of the smooth connected components $W_i$ of the algebraic hypersurface $Y_0(P)$.
\et
\pr
In order to avoid complicated expressions, we use below constants $C_q(d,n)$, depending only on $d$ and $n$, not specifying some of them explicitly.

\medskip

The following lemma provides a bound for the norms of the intermediate derivatives $M_q(f)$ of a smooth $f$ through $M_0(f)$ and $M_{d+1}(f)$.

\bl\label{lem:norms.der}
Let $f$ be a $C^{d+1}$-smooth function on $B^n.$ Then for $k=1,2,\ldots,d$ we have
$$
M_k(f)\le C_1(n,d)M_0(f)+ C_2(n,d)M_{d+1}(f).
$$
\el
\pr
Let $P=P_d(f)$ be the Taylor polynomial at the origin of degree $d$ of $f$. By Taylor's formula we have for $x \in B^n$
$$
|f(x)-P(x)|\le \frac{1}{(d+1)!}M_{d+1}(f).
$$
We conclude that $M_0(P)\le M_0(f)+\frac{1}{(d+1)!}M_{d+1}(f)$. Next we use the equivalence of all the norms on the finite-dimensional space of polynomials of degree $d$, and obtain, for $k=1,2,\ldots,d,$
$$
M_k(P)\le \bar C(n,d)M_0(P)\le \bar C(n,d)[M_0(f)+\frac{1}{(d+1)!}M_{d+1}(f)].
$$
Finally, once more using Taylor's formula, we get for $k=1,2,\ldots,d,$
$$
M_k(f)\le  M_k(P)+\frac{1}{(d-k+1)!} M_{d+1}(f)\le  C_1(n,d)M_0(f)+ C_2(n,d)M_{d+1}(f),
$$
where $C_1(n,d)=\bar C(n,d), \ \ C_2(n,d)= \bar C(n,d)\frac{1}{(d+1)!}+\frac{1}{(d-k+1)!}.$ This completes the proof of Lemma \ref{lem:norms.der}. $\square$

\smallskip

In particular, under our initial assumptions that $M_0(f)=1$ and $M_{d+1}(f)\le T \le 1$ we have
$$
M_2(f)\le C_1(n,d)+C_2(n,d):=C_3(n,d).
$$
Put $\delta=\frac{\gamma}{3C_3}$. Then, in a $\delta$-neighborhood $U_\delta$ of $Y_0(f)$ we have $||grad \ f(x)||\ge \frac{1}{2}\gamma$.

\smallskip

Next we consider a vectorfield $v(x)$ in $U_\delta$ defined by $v(x)=\frac{grad \ f(x)}{||grad \ f(x)||^2}$. For $x\in U_\delta$ we have $||v(x)||\le \frac{2}{\gamma}$.

\medskip

The derivative of $f$ in the direction of $v(x)$ satisfies the identity
$$
\frac{df}{dv}(x)=<v(x), grad \ f(x)>=\frac{<grad \ f(x), grad \ f(x)>}{||grad \ f(x)||^2}=1.
$$
For each $y\in Y_0(f)$ denote by $\zeta(y,t)$ the trajectory of the vectorfield $v(x)$, satisfying $\zeta(y,0)=y.$ Since for $x\in U_\delta$ we have $||v(x)||\le \frac{2}{\gamma}$, for each $t$ with $|t|\le \frac{\delta \gamma}{2}=\frac{\gamma^2}{2C_3}:=\eta$ the trajectory $\zeta(y,t)$ remains in $U_\delta$, and hence it is well-defined.

\medskip

Now we consider the ``normal bundle'' mapping $\Psi: G\to U_\delta$ of the product $G=Y_0(f)\times [-\eta,\eta]$ into $U_\delta$, defined by
$$
\Psi(y,t)=\zeta(y,t).
$$
By the construction, $\Psi$ satisfies $f(\Psi(y,t))=t.$ By the uniqueness and dependence on the initial data of the trajectories $\zeta(y,t),$ the mapping $\Psi$ provides a diffeomorphism of $G$ with its image, which is the level strip $Q_\eta=\{x\in U_\delta, \ |f(x)|\le \eta\}$.

\medskip

Let us return now to the Taylor polynomial $P$ of degree $d$ of $f$ at the origin. By the remainder formula, and since we assume that $M_{d+1}(f)\le T$,  we have $M_0(f-P)\le \frac{T}{(d+1)!} \le \frac{1}{2}\eta$. Also by the Taylor formula we have $M_1(f-P)\le \frac{T}{d!} \le \frac{\gamma^2}{4C_3}$. We conclude that for $t\in [-\eta,\eta],$ and for each $y\in Y_0(f),$ along the trajectory $\zeta(y,t)$ we have $|\frac{dP}{dt}-1|\le \frac{1}{2}$.

\medskip

Therefore, for each $y\in Y_0(f),$ along the trajectory $\zeta(y,t)$ the polynomial $P(\zeta(y,t))$ has exactly one simple zero at a certain $t(y)\in [-\eta,\eta].$ By the implicit function theorem, the function $t(y)$ is $C^{d}$-smooth.

\medskip

Summarizing, we conclude that the part $\O$ of the zero set $Y_0(P)$ of $P$, which is contained in the level strip $Q_\eta$, is given in the coordinates $y,t$ on $Q_\eta$ as the graph of the smooth function $t(y)$ on $Y_0(f)$. Hence
$\O=\cup_{i=1}^q \O_i ,$ where each $\O_i$ is the graph $t(y)$ on $V_i$. We immediately conclude also that $\O_i$ is diffeomorphic to $V_i, \ i=1,\ldots,q$. But in fact, the formula $t_\tau(y)=\tau t(y), \ \tau\in [0,1],$ provides a smooth isotopy between $\O$ and $Y_0(f)$. This completes the proof of Theorem \ref{thm:rigidity.topology}. $\square$

\medskip

The corresponding rigidity statement is

\bt\label{thm:rigid.topol}
Let $f$ be as above, with $0$ a $\gamma$-regular value of $f$. If $Y_0(f)$ is not smoothly isotopic to a certain union $W=\cup_{i=1}^q W_i $ of smooth connected components $W_i$ of an algebraic hypersurface $Y_0(P)$, with $P$ a polynomial of degree $d$, then
$$
M_{d+1}(f)\ge T=T(\gamma)=\min \{1,\frac{d!\gamma^2}{4C_3}\}.
$$
\et

\medskip

\noindent{\bf Remark} ``No isotopy'' condition of Theorem \ref{thm:rigid.topol} can be weakened in many forms, in particular, in terms of the Betti numbers of the components of $Y_0(f)$, of their mutual position, etc.


\bibliographystyle{amsplain}

\begin{thebibliography}{10}

\bibitem{Bih} F. Bihan, Asymptotic behavior of Betti numbers of real algebraic surfaces,
Comment. Math. Helv. 78 (2003), no. 2, 227--244.

\bibitem{Bru.Yom} A. Brudnyi, Y. Yomdin, Norming sets, and related Remez-type inequalities,
J. Aust. Math. Soc. 100 (2016) 163--181

\bibitem{Bru.Gan} Y. Brudnyi, M. Ganzburg, On an extremal problem
for polynomials of $n$ variables, {\sl Math. USSR Izv.} {\bf 37}
(1973), 344-355.

\bibitem{Bru.Shv} Y. Brudnyi and P. Shvartsman,  Generalizations of Whitney’s extension theorem, IMRN 1994, no. 3, 129--139.


\bibitem{Erd} T. Erdelyi, Remez-type inequalities and their
applications, {\sl J. Comp. Appl. Math.} {\bf 47} (1993) 167-209.

\bibitem{Fef} C. Fefferman, Whitney’s extension problem for $C^m$, Ann. of Math. 164 (2006), 313--359.

\bibitem{Fef.Kla} C. Fefferman, B. Klartag, Fitting a $C^m$-smooth function to
data. Part I, {\sl Ann. of Math.} (2) 169 (2009), 315–346. Part II, {\sl Rev. Mat.
Iberoam.} 25 (2009), 49–273.

\bibitem{Ler.Ste} A. Lerario, M. Stecconi, What is the degree of a smooth hypersurface?
arXiv:2010.14553v1.


\bibitem{Rem} E. J. Remez, Sur une propriete des polynomes de
Tchebycheff, {\sl Comm. Inst. Sci. Kharkov} {\bf 13} (1936) 93-95.

\bibitem{Whi1} H. Whitney, Analytic extensions of differentiable functions defined in closed sets, Trans.
Amer. Math. Soc. 36 (1934), 63--89.

\bibitem{Whi2} H. Whitney, Differentiable functions defined in closed sets. I, Trans. Amer. Math. Soc. 36
(1934), 369--387.

\bibitem{Whi3} H. Whitney, Functions differentiable on the boundaries of regions, Ann. of Math. 35 (1934),
482--485.

\bibitem{Yom1} Y. Yomdin, The set of zeroes of an ``almost polynomial'' function,
{\sl Proc. AMS}, Vol. 90, No. 4 (1984), 538-542.

\bibitem{Yom.Com} Y. Yomdin, G. Comte, Tame geometry with application in smooth analysis,
Springer, LNM 1843, 2004.

\bibitem{Yom} Y. Yomdin, Remez-Type Inequality for Discrete Sets,
{\sl Isr. J. of Math.}, Vol 186 (Nov. 2011), 45-60.

\bibitem{Yom2} Y. Yomdin, Remez-Type Inequality for smooth functions, in
``Geometry and its applications'', V. Rovenski, ‎P. Walczak, Editors, Springer, Cham. 2014, 235-243.

\bibitem{Yom3} Y. Yomdin, Smooth rigidity and Remez-Type Inequalities, Anal.Math.Phys. 11, 89 (2021).

\bibitem{Yom4} Y. Yomdin, Smooth rigidity via probe curves, in preparation.







\end{thebibliography}

\end{document}